# SEVERAL CLASSICAL IDENTITIES VIA MELLIN'S TRANSFORM

## KHRISTO N. BOYADZHIEV

Department of mathematics, Ohio Northern University
Ada, Ohio, 45810, USA

E-mail: k-boyadzhiev@onu.edu

**DOI:**

**Summary.** We present a summation rule using Mellin's transform to give short proofs of some important classical relations between special functions and Bernoulli and Euler polynomials. For example, the values of the Hurwitz zeta function at the negative integers are expressed in terms of Bernoulli polynomials. We also show identities involving exponential and Hermite polynomials.

## 1. INTRODUCTION

Throughout we use the notation $(a) = a + it, t \in \mathbb{R}$ for the vertical line with abscissa $0 < a < 1$, oriented from minus to plus infinity. First we recall the formulas for the Mellin transform

$$G(s) = \int_0^\infty x^{s-1} g(x) dx$$

and its inverse

$$g(x) = \frac{1}{2\pi i} \int_{(a)} x^{-s} G(s) ds, \ x > 0. \qquad (1)$$

**Definition**. For $n > 0$ consider integrals of the form

$$\int_{L(n)} x^{-s} G(s) ds, \ x > 0$$

where $L(n)$ consists of the line segment $[-ni, ni]$ together with the semicircle $R(n)$ in the left half plane for which the line segment is the diagonal. If the integrals

$$\int_{R(n)} x^{-s} G(s) ds$$

approach zero when $n \to \infty$, we say that the line of integration in (1) can be closed to the left. In a similar manner we define integrals where the line of integration can be closed to the right.

**Proposition.** Suppose that the function $G(s)$ is meromorphic on the half plane $Re(s) < a + \varepsilon$ for some small $\varepsilon > 0$, and has only simple poles at $s = 0, -1, -2, \ldots$, with residues $c_0, c_1, c_2, \ldots$. If the line of integration in (1) can be closed to the left, the residue theorem provides the representation

$$g(x) = \sum_{n=0}^{\infty} c_n x^n \tag{2}$$

for the function $g(x)$ from (1), i.e. this function is a power series. If now $f(s)$ is an appropriate holomorphic function on $Re(s) < a + \varepsilon$ without poles, we can write

$$\frac{1}{2\pi i} \int_{(a)} x^{-s} f(s) G(s) \, ds = \sum_{n=0}^{\infty} c_n f(-n) x^n, \tag{3}$$

when the power series on the right side converges.

Formulas of this type are used for summation of series or interpolation, and are present in many publications (see [2, 4, 9, 10] and the references there).

In this note we focus on a special area of applications for the proposition - obtaining some classical identities by using Mellin inversion and comparing coefficients. In order to keep the paper short we omit details and do not discuss convergence of some integrals and series. The validity of such formulas is considered in [4, 9, 10].

The illustration of the method is given in the following examples.

## 2  EXAMPLES

Remind that the residues of the gamma function at zero and the negative integers are given by $\operatorname{Res}(\Gamma, -n) = \dfrac{(-1)^n}{n!}$ for $n = 0, 1, 2, \ldots$. Also $\Gamma(s)$ has rapid decay on vertical lines. We have the estimate ([15, (20)]

$$|\Gamma(a+it)| \sim \sqrt{2\pi} \, |t|^{a-\frac{1}{2}} e^{-\frac{\pi}{2}|t|}, a \in \mathbb{R}$$

when $|t| \to \infty$. The estimate helps for the convergence of our integrals.

The above facts will be used in the following examples. We note that in these examples the lines of integration can be closed to the left (proofs are standard, using the growth estimate for $\Gamma(s)$). Note also that when we replace $G(s)$ by $\Gamma(s)$ in (3) we have

$$\frac{1}{2\pi i} \int_{(a)} x^{-s} f(s) \Gamma(s) \, ds = \sum_{n=0}^{\infty} \frac{(-1)^n}{n!} f(-n) x^n. \tag{4}$$

**Example 1.** Ramanujan's Master Theorem. As Hardy writes in [8], Ramanujan was very fond of his integral formula [8, p.186]

$$\int_0^{\infty} x^{s-1} \left\{ \sum_{n=0}^{\infty} \frac{f(n)(-x)^n}{n!} \right\} dx = f(-s) \Gamma(s) \tag{5}$$

and used it for many applications. Berndt rightly calls it Ramanujan's Master Theorem [2, Entry 11, p.105]. Details, comments, and applications of (5) are given in these two books and also in [1, 4, 6]. Clearly, after replacing $f(s)$ by $f(-s)$ equation (5) turns into (4) after Mellin inversion. Ramanujan did not use the residue theorem for his proof but only standard calculus (see [2, p.106]).

**Example 2.** The Hurwitz zeta function is defined by

$$\zeta(s,z) = \sum_{n=0}^{\infty} \frac{1}{(n+z)^s} \quad (\mathrm{Re}(z) > 0, \mathrm{Re}(s) > 1)$$

with integral representation

$$\zeta(s,z)\Gamma(s) = \int_0^{\infty} t^{s-1} \frac{e^{t(1-z)}}{e^t - 1} dt \quad (\mathrm{Re}(s) > 1)$$

When $z = 1$, $\zeta(s,1) = \zeta(s)$ is the Riemann zeta function [5].

We have also the modified integral representation (argument is the same as on pp. 61-62 in [13])

$$\zeta(s,z)\Gamma(s) = \int_0^{\infty} t^{s-1} \left( \frac{e^{t(1-z)}}{e^t - 1} - \frac{1}{t} \right) dt \quad (0 < \mathrm{Re}(s) < 1)$$

which is a Mellin transform formula. By Mellin inversion

$$\frac{e^{t(1-z)}}{e^t - 1} - \frac{1}{t} = \frac{1}{2\pi i} \int_{(a)} t^{-s} \zeta(s,z) \Gamma(s) ds = \sum_{n=0}^{\infty} \frac{(-1)^n}{n!} \zeta(-n, z)$$

At the same time, the Bernoulli polynomials $B_n(z)$ have the generating function

$$\frac{x e^{xz}}{e^x - 1} = \sum_{k=0}^{\infty} \frac{B_k(z)}{k!} x^k$$

and from this

$$g(x,z) \equiv \frac{e^{x(1-z)}}{e^x - 1} - \frac{1}{x} = \sum_{k=1}^{\infty} \frac{B_k(1-z)}{k!} x^{k-1} = \sum_{n=0}^{\infty} \frac{B_{n+1}(1-z)}{(n+1)!} x^n.$$

Therefore, by comparing coefficients for $n = 0, 1, \ldots$ we find the classical formula

$$\zeta(-n, z) = \frac{(-1)^n B_{n+1}(1-z)}{n+1} = -\frac{B_{n+1}(z)}{n+1}$$

(using the property $(-1)^n B_n(1-z) = B_n(z)$, so that $(-1)^n B_{n+1}(1-z) = -B_{n+1}(z)$).

In particular,

$$\zeta(0, z) = -B_1(z) = \frac{1}{2} - z.$$

For the Bernoulli numbers $B_n = B_n(0) = (-1)^n B_n(1)$ we have

$$\zeta(-n) = \zeta(-n, 1) = \frac{(-1)^n B_{n+1}(0)}{n+1} = \frac{(-1)^n B_{n+1}}{n+1}.$$

The odd Bernoulli numbers are zeros except $B_1 = -1/2$. Thus $\zeta(0) = -1/2$.

Next we give a new proof of a result of Kenneth Williams and Zhang Nan-Yue [14].

**Example 3.** Consider now the alternating Hurwitz zeta function

$$\eta(s,z) = \sum_{n=0}^{\infty} \frac{(-1)^n}{(n+z)^s} \quad (\operatorname{Re}(z) > 0, \operatorname{Re}(s) > 0)$$

which extends to the entire complex plane as analytic in the variable $s$ and has the integral representation

$$\eta(s,z)\Gamma(s) = \int_0^{\infty} t^{s-1} \frac{e^{t(1-z)}}{e^t + 1} dt \quad (\operatorname{Re}(s) > 0).$$

By inversion

$$\frac{e^{t(1-z)}}{e^t + 1} = \frac{1}{2\pi i} \int_{(a)} t^{-s} \eta(s,z)\Gamma(s) ds = \sum_{n=0}^{\infty} \frac{(-1)^n}{n!} \eta(-n,z) x^n.$$

Euler's polynomials $E_n(z)$ are defined by the generating function

$$\frac{2e^{tx}}{e^t + 1} = \sum_{n=0}^{\infty} E_n(x) \frac{t^n}{n!} \quad (|t| < \pi).$$

Comparing coefficients gives

$$\eta(-n,z) = \frac{(-1)^n}{2} E_n(1-z) = \frac{1}{2} E_n(z)$$

by the property $E_n(1-z) = (-1)^n E_n(z)$.

**Example 4.** Euler worked with the function

$$L(s) = \sum_{n=0}^{\infty} \frac{(-1)^n}{(2n+1)^s} \quad (\operatorname{Re} s > 0)$$

which we call here Euler's $L$-function (sometimes it is called Dirichlet's $L$-function). This function has the integral representation

$$2\Gamma(s)L(s) = \int_0^{\infty} \frac{x^{s-1}}{\cosh x} dx.$$

It also has analytic extension on the complex plane.

By Mellin inversion and using equation (4)

$$\frac{1}{2\cosh(x)} = \frac{1}{2\pi i} \int_{(a)} x^{-s} L(s)\Gamma(s) ds = \sum_{n=0}^{\infty} \frac{(-1)^n}{n!} L(-n) x^n$$

Euler's numbers $E_n$ are defined by the generating function

$$\frac{1}{\cosh x} = \sum_{n=0}^{\infty} \frac{E_n}{n!} x^n.$$

This function is even, so the Euler numbers with odd indices are zeros. By comparing coefficients we find

$$L(-2n) = \frac{1}{2} E_{2n} \quad (n = 0, 1, 2, ...).$$

**Example 5.** The exponential polynomials $\varphi_n$ are defined by the generating function

$$e^{z(e^x - 1)} = \sum_{n=0}^{\infty} \varphi_n(z) \frac{x^n}{n!}$$

(see [3, 12]). We will use the function

$$\psi(x, z) = e^{z(e^{-x} - 1)} = \sum_{n=0}^{\infty} (-1)^n \varphi_n(z) \frac{x^n}{n!}$$

which has Mellin transform

$$\Psi(s, z) = \int_0^{\infty} x^{s-1} \psi(x, z) \, dx = e^{-z} \int_0^{\infty} x^{s-1} e^{z e^{-x}} \, dx$$

$$= e^{-z} \sum_{n=0}^{\infty} \frac{z^n}{n!} \left\{ \int_0^{\infty} x^{s-1} e^{-nx} dx \right\} = e^{-z} \Gamma(s) \sum_{n=0}^{\infty} \frac{z^n}{n! n^s} = \Gamma(s) f(s, z)$$

with

$$f(s, z) = e^{-z} \sum_{n=0}^{\infty} \frac{z^n}{n! n^s}.$$

Now we have from equation (4)

$$\psi(x, z) = \frac{1}{2\pi i} \int_{(a)} x^{-s} \Psi(x, s) ds = \frac{1}{2\pi i} \int_{(a)} x^{-s} f(s, z) \Gamma(s) ds = \sum_{n=0}^{\infty} \frac{(-1)^n}{n!} f(-n, z) x^n$$

Thus for $n \geq 0$ (with the agreent $0^0 = 1$)

$$\varphi_n(z) = f(-n, z) = e^{-z} \sum_{k=0}^{\infty} \frac{k^n}{k!} z^k \tag{6}$$

which is one of the fundamental properties of the exponential polynomials.

Note that identity (6) can be used for a meaningful extension of $\varphi_n(z)$ to $\varphi_\lambda(z)$ with non-integer index $\lambda$. For instance, we have (with $0^\lambda = 0$)

$$\varphi_\lambda(z) = e^{-z} \sum_{k=1}^{\infty} \frac{k^\lambda}{k!} z^k.$$

**Example 6**. The generating function for the Hermite polynomials is

$$\psi(x,z) = e^{2xz-x^2} = \sum_{n=0}^{\infty} H_n(z)\frac{x^n}{n!}$$

with Mellin transform [11, p. 27]

$$\Psi(s,z) = e^{\frac{z^2}{2}} 2^{-s/2} D_{-s}(\sqrt{2}z)\Gamma(s).$$

where $D_p$ are the parabolic cylinder functions [7, pp. 1065-1067]. Let now

$$f(s,z) = e^{\frac{z^2}{2}} 2^{-s/2} D_{-s}(\sqrt{2}z).$$

Then $H_n(z) = (-1)^n f(-n,z)$ and since $H_n(z) = (-1)^n H_n(-z)$ one finds the classical result

$$H_n(z) = 2^{\frac{n}{2}} e^{\frac{z^2}{2}} D_n(\sqrt{2}z)$$

([7, entry 9.253, p.1067]).

## 6. CONCLUSIONS

In this note we presented a rule how Mellin's transform can be used to give short proofs of several classical identities connecting, in particular, the Bernoulli and Euler polynomials to the values of the Hurwitz and alternating Hurwitz functions at the negative integers. We also proved identities for the exponential and Hermite polynomials.